\documentclass[12pt]{elsarticle}
\usepackage{times}
\usepackage{lineno,hyperref}
\modulolinenumbers[5]

\usepackage{color}
\usepackage[noend]{algpseudocode}
\usepackage{algorithm}
\usepackage{amsmath,bm}
\usepackage{amssymb}
\usepackage{array}
\usepackage{geometry}
\usepackage{setspace}
\journal{REIRINA Summer Internship PMU Group}

\begin{document}

\begin{frontmatter}

\title{Bayesian Estimation Based Load Modeling Report}

\author{Chang Fu}
\ead{chang.fu@wayne.edu}

\author{Zhe Yu}
\ead{zhe.yu@geirina.net}

\author{Di Shi}
\ead{di.shi@geirina.net}

\address{GEIRI North America, San Jose, CA, 95134, USA}

%
%
%
%
%
%





\end{frontmatter}


\section{Background and Motivations}
\subsection{Frequentist vs Bayesian probability}
Classical probability and Bayesian probability are the two major interpretations of the probability theory. The classical probability defines an event's probability as the limit of its relative frequency $k$ in a large number of trials $n$, i.e.
\begin{equation}
p=\lim_{n \to \infty}\frac{k}{n}
\end{equation}\par
The most frequently used classical parameter estimation approach is maximum likelihood (ML) method, where it finds the parameters values that maximize the likelihood function for the observed data and then estimates the errors of these estimates \cite{book1}.  In classical probability, the parameters are not associated with probability distributions. \par
In contrast, Bayesian probability belongs to the category of evidential probabilities. To evaluate the probability of a hypothesis,  Bayesian probability specifies some prior probability, which is then updated to a posterior probability in the light of new, relevant data (likelihood) \cite{07}. The formula of Bayesian estimation can be written as:
\begin{eqnarray}
p(\theta|x_i)=\frac{p(\theta)\cdot p(x_i|\theta)}{\sum p(\theta)\cdot p(x_i|\theta)}
\end{eqnarray}
where $p(\theta)$ is the prior, $p(\theta|x_i)$ is the posterior. $p(x_i|\theta)$ is the likelihood. The denominator is called normalize factor. Bayesian estimation gives the probability density function (PDF) of a parameter which is given by some already known observations. It is the product of the  combination of the prior information as well as the observations of happened events (likelihood). Different with a specific number given by ML, a PDF is able to provide more information and can be used in forecasting, optimization, etc..
\subsection{Composite Load Modeling Methods}
Load modeling is a very traditional topic in power system and has been widely studied in \cite{08,09,gibbsmjin, expzip1,10,30,31,32} in  past couple of decades. An accurate load model has been shown of great importance in voltage stability studies \cite{book2} as well as the planning, operation and control of power systems \cite{11}.The emerging of smart grid technologies such as distributed generations (DGs) \cite{12}, load/demand side management, optimal dispatch \cite{13}, etc,  make it more and more desired for establishing an appropriate load model in order to achieve a more accurate and reliable system.\par
\begin{figure}[!h]
\centering
    \includegraphics[width=0.8\linewidth,]{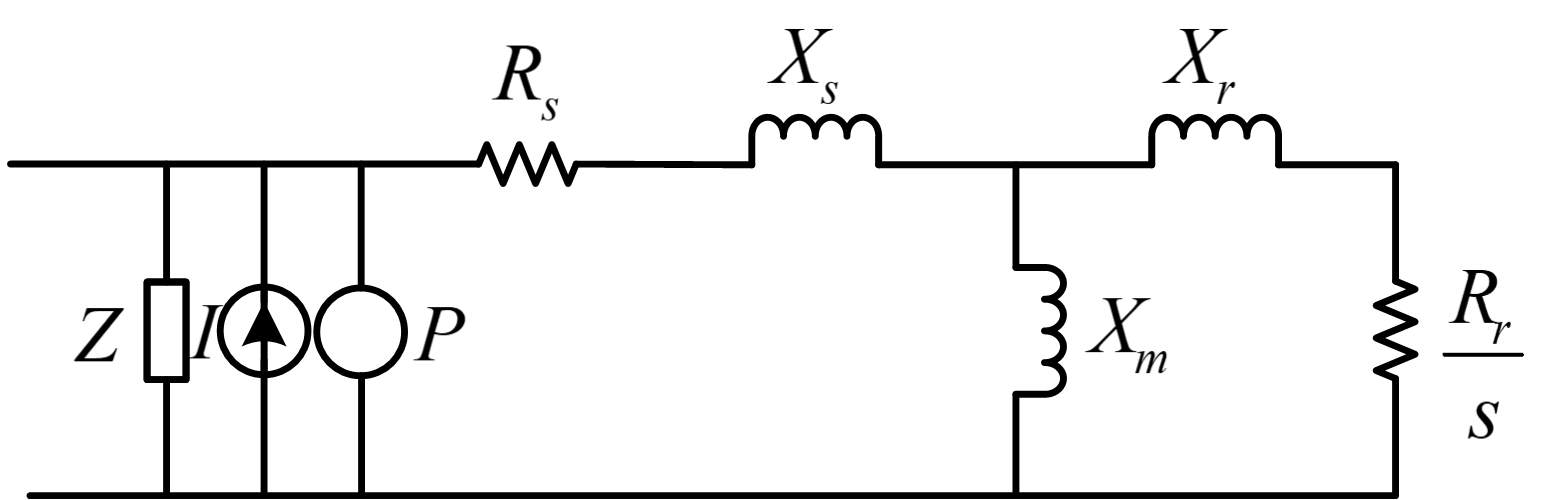}
\caption{The composite model of ZIP+IM.}
\label{composite}
\end{figure}

Based on different mathematical presentations, the load models can be categorized into static model and dynamic model. Static models include ZIP model, exponential model, frequency dependent model, etc.. Examples of widely used dynamic models including induction motor (IM), exponential recovery load model (ERL) \cite{08}. In most cases, researches were focused on  composite load modeling consisting  static and dynamic load models \cite{08}. It is reported in \cite{08} that ZIP+IM model can be used with various conditions, locations and compositions. A simplified composite model using ZIP+IM model can be found in Fig. \ref{composite}. The task reported here is how to estimate/identify parameters for ZIP+IM model.\par
Component based and measurement based approaches are the two main load model parameter identification methods which have been widely used nowadays \cite{16,17}. Component based approach requires information of load compositions. As introduction in \cite{08}, this approach requires models of individual components, component compositions (consumption of different load type in percentage), as well as class composition (residential, industrial or commercial load percentage).  Measurement-based approaches calculate the parameters using algorithms such as least-squares (LS), genetic algorithm (GA) based on selected model structure by minimize the objective function:
\begin{equation}
\textup{min}\ \frac{1}{n}\sum_{i=1}^n[(P_i^m-P_i^e)^2+(Q_i^m-Q_i^e)^2]
\end{equation}\par
where superscript "$m$" indicates the measured value, "$e$" presents the estimated value, respectively. The measurement should be obtained in different conditions and disturbances. Different techniques have been used in \cite{16,20,21} to identify the parameters for composite model. Other methods such as Hammerstein model based real-time system identification for dynamic load was proposed in \cite{18}. An Artificial neural network (ANN)-based modeling were proposed in \cite{19}.
\subsection{Motivations for using Bayesian Estimation}
Traditional parameters estimation methods, e.g. measurement based approach, can only be obtained at specific locations within certain period of time. Thus, in \cite{15}, robust time-varying load parameters identification approach was proposed by using  batch-model regression in order to obtain the updated system parameters. Moreover, measurement based approach is highly dependent on measurements from the meters. Measurement anomalies may effect the robustness of the estimations in both time-varying and time-independent load models. Consequences of measurement noises is another significant issue need to be addressed when implementing corresponding methods \cite{14}.  
Moreover, component based approach highly relies on the characteristics of individual load components and compositions of  load . In \cite{expzip1}, ZIP coefficients for the most commonly used electrical appliances were determined in order to get an accurate model before there critical values can be applied into the model.  The process of determination of different components involves characterization of the consumers, clustering, etc.. Since the data are highly dependent on the type of the day, weather, service class, it is a very computationally wasted approach to get ZIP coefficients for these appliances. Furthermore, it is introduced in \cite{08}, the accuracy of such approach is still insufficient. Difficulties involved in obtaining the comprehensive load composition information is another fact need to be carefully considered when implementing this method. Therefore, some efforts were made and a Bayesian based method was carried out in this report.\par 
Bayesian estimation (BE) based composite load parameter identification approach can successfully overcome the aforementioned disadvantages in measurement-based and component-based approaches. First, BE is a time-varying estimation method. Parameters can be updated according to the requirements. Furthermore, BE does not require the complicated procedure used in component based approach. The information of the composition of the loads and the coefficients of the appliances are not necessarily needed. Furthermore, BE is a robust estimation method due to its statistic characteristics: measurement anomalies will not significantly effect the results if the number of samples is large enough, and the measurement error will not effect the estimation results since in each step of sampling the estimated parameter is drawn from a normal distribution calculated according to the data and eventually formed a distribution of the target parameter. 

\section{ZIP Model}

The most commonly used composite model in industry is ZIP+IM model.  ZIP model represents the static part of a composite load model. Detailed format of a ZIP model can be found in \ref{zipp} and \ref{zipq}:
\begin{eqnarray}
\label{zipp}
P_{ZIP}&=&P_0(\alpha_1\bar{V}^2+\alpha_2\bar{V}+\alpha_3)\\
\label{zipq}
Q_{ZIP}&=&Q_0(\alpha_4\bar{V}^2+\alpha_5\bar{V}+\alpha_6)\\
\sum_{i=1}^3{\alpha_i}&=&\sum_{i=4}^6{\alpha_i}=1
\end{eqnarray} 
where $\bar{V}=\dfrac{V}{V_0}$, $V_0$ is the voltage at the initial value, $V$ is the measured value. Coefficients of ZIP model were fitted to the measured data with the constrains that the sum of all coefficients is 1.  ZIP coefficients for reactive power are typically large because of the high load power factors and high partial derivative \cite{zipcoefficient}.

\subsection{Markov Chain Monte Carlo (MCMC) Method}
\begin{figure}[!h]
\centering
    \includegraphics[width=0.55\linewidth,]{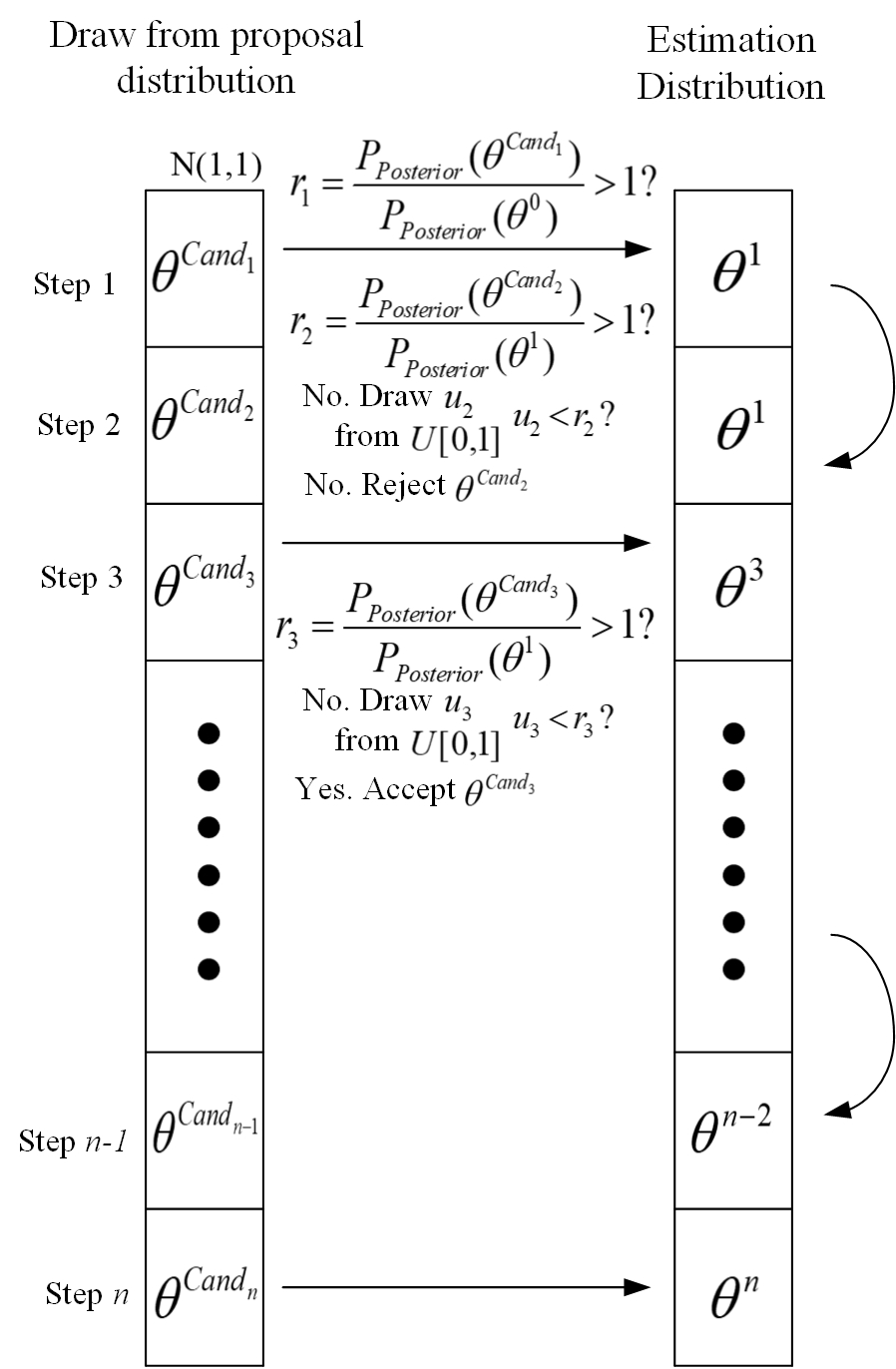}
\caption{MCMC introduction.}
\label{mcmcintro}
\end{figure}

One of the Bayesian estimation based parameter estimation approach is MCMC method. It is a method that repeatedly draws random values for the parameters of a distribution based on the current values. The reason that we can use Bayesian estimation based methods is because Markov Chain does not depends on the initial status. It will eventually converge to its final state. We using Markov Chain to approach the real value of the parameter and finally gives the estimation distribution of this parameter. Detailed steps of implementing a MCMC is shown in Fig. \ref{mcmcintro}. This specific method is called Metropolis-Hastings (MH) method, it is a updated version of MCMC since it accelerate the acceptance in this method. The burning-in period in MH is much shorter than that in MCMC.\par
  There are some assumptions were made before MH was used: (1) The active/reactive power ($P, Q$), initial active/reactive power ($P_0,Q_0$), bus voltage ($V$), initial bus voltage ($V_0$) are all known. (2) The measurements of aforementioned values were independent identical distributed (i.i.d.). (3) The distributions of \ref{zipp} and \ref{zipq} follow normal distribution.\par
MH method works well with only one parameter unknown, and our specific problem can be described as follows:
 For a $ZIP$ model with the form $P=P_0(\alpha_1\bar{V}^2+\alpha_2\bar{V}+\alpha_3)$, if we only consider there is one type of load (e.g., only constant I), the equation can be written as $P=P_0(\alpha_2\bar{V})$. Now we want to use the MCMC method to estimate the probability distribution of $\alpha_2$.\par
The steps are as follows:\par
1. generate a proposal distribution, such as a normal distribution :
\begin{equation}
 X \sim \mathcal{N}(\mu,\,\sigma^{2})
  \end{equation}
with $\mu$ as the mean, $\sigma$ as the variance. We can chose any arbitrary parameters for these parameters. In this case, just use 
\begin{equation}
X \sim \mathcal{N}(1,1)
\label{normone}
\end{equation}
 for instance.\par
2. Draw a random number $\theta_1$ from \ref{normone}. Let the initial condition as $\theta_0=0.1$ (this value is arbitrary), calculate the ratio $r$, where
\begin{equation}
r=\frac{P_{posterior}(\theta_t)}{P_{posterior}(\theta_{t-1})}=\frac{P_{prior}(\theta_t)\cdot P_{likelihood}(\theta_t)}{P_{prior}(\theta_{t-1})\cdot P_{likelihood}(\theta_{t-1})}
\label{calculater}
\end{equation}
In \ref{calculater}, the prior is an approximate distribution given by the expert experience, (e.g., it can be a normal distribution with $\mu_p$ and $\sigma_p$, which is $X \sim \mathcal{N}(\mu_p,\,\sigma_p^{2})$). The likelihood is the observation of  an event happens given by the parameter $\theta$.\par
More specifically, in the first step,
\begin{equation}
r_1=\frac{P_{posterior}(\theta_1)}{P_{posterior}(\theta_{0})}=\frac{P_{prior}(\theta_1)\cdot P_{likelihood}(\theta_1)}{P_{prior}(\theta_{0})\cdot P_{likelihood}(\theta_{0})}
\end{equation}
If this value is larger than 1, $\theta_1$ is accepted. If the value is smaller than 1, a random number $u_1$ will be chosen from a uniform distribution $U \sim \mathcal{U}[0,1]$, then $r_1$ and $u_1$ are compared. Similarly, if $r_1>u_1$, $\theta_1$ will be accepted, otherwise reject $\theta_1$. Therefore, in next step, the new parameter used in the denominator will be update to $\theta_1$ if accepted, or remain $\theta_0$ if rejected. Meanwhile, a new $\theta_2$ will be drawn from \ref{normone}. $r_2$ will be calculated by using the following equation:
\begin{equation}
r_2=\frac{P_{posterior}(\theta_2)}{P_{posterior}(\theta_{1})}=\frac{P_{prior}(\theta_2)\cdot P_{likelihood}(\theta_2)}{P_{prior}(\theta_{1})\cdot P_{likelihood}(\theta_{1})}
\end{equation}
In \ref{calculater}, $P_{prior}(\theta_t)$ can be calculated from prior distribution using $\theta_t$ (this value is drawn from the proposal distribution). $P_{likelihood}(\theta_t)$ is calculated by using the number of the event happens in all observation that equals $\theta_t$ ($\theta_t$ is not a probability, it is just a possible number in the norm distribution) divided by the total number of the observation.
\makeatletter
\def\algbackskip{\hskip-\ALG@thistlm}
\makeatother

\begin{algorithm}[!h]
\caption{Metropolis-Hastings Algorithm}\label{mcmc}
\begin{algorithmic}[1]
\State $\text{Initialize}\ \theta^{(0)} \ \text{from proposal distribution}\ p(\theta)$ 
\For {$\text{iteration}\ i=1,2,...$}
\State $\text{Propose}\  \theta^{Cand}\sim p(\theta^{(i)}|\theta^{(i-1)})$
\State $\text{Calculate}\  r_i=\dfrac{P_{Posterior}(\theta^{Cand})}{P_{Posterior}(\theta^{(i-1)})}$
\If {$r_i>1$}
\State $\text{Accept the proposal}\ \theta^{(i)}\ \gets \theta^{Cand}$
\Else \ {$\text{Draw a random number}\ u^{(i)}\ \text{from}\ \mathcal{U}[0,1]$}
\If {$u^{(i)}>r^{(i)}$}
\State $\text{Reject the proposal}\ \theta^{(i)}\ \gets \theta^{(i-1)}$
\Else \ {$\text{Accept the proposal}\ \theta^{(i)}\ \gets \theta^{Cand}$}
\EndIf
\EndIf
\EndFor
\end{algorithmic}
\end{algorithm}

\subsection{Gibbs Sampling (GS)}
 In most cases, the three parameters of (\ref{zipp}) are all unknown. Though the MH method can be extended to 3 parameters, but finding an appropriate acceptance rate is a big issue in implementing MH with multi-variate \cite{book1}. Therefore, Gibbs sampling will be used for multi-variate parameters estimations. The principle in Gibbs sampling is that only one parameter will be sampled/estimated with remaining parameters fixed in each step. In next step, the parameter estimated in previous step will be used to replace the corresponding one in current step and a new parameter will be sampled. The pseudo code can be found as follows:
\makeatletter
\def\algbackskip{\hskip-\ALG@thistlm}
\makeatother

\begin{algorithm}[!h]
\caption{Gibbs Sampling}\label{gibbs}
\begin{algorithmic}[1]
\State $\text{Initialize}\ x^{(0)} \sim q(x)$ 
\For {$\text{iteration}\ i=1,2,...$}
\State $x_1^{(i)}\sim p(x_1|x_2^{(i-1)},x_3^{(i-1)},...,x_n^{(i-1)})$
\State $x_2^{(i)}\sim p(x_2|x_1^{(i)},x_3^{(i-1)},...,x_n^{(i-1)})$
\State $\vdots$
\State $x_n^{(i)}\sim p(x_n|x_1^{(i)},x_3^{(i)},...,x_{n-1}^{(i)})$
\EndFor
\end{algorithmic}
\end{algorithm}
Before implementing Gibbs Sampling method, define:\par
1. $y_i=\dfrac{P_i}{P_0}$, $x_i=\dfrac{V_i}{V_0}$. Where $P_0$ is the initial active power without connecting the ZIP load, $P_i$ is the active power with the target load connected in the $i$th experiment. Their corresponding voltage are defined as $V_0$ and $V_i$, respectively. \par
2. $y_i\sim \mathcal{N}(\alpha_1x_i^2+\alpha_2x_i+\alpha_3,1/\tau)$. Since $\alpha_1+\alpha_2+\alpha_3=1$, the number of unknown variable is reduced to 2.\par
3. Total number of $n$ experiments were made, which means there are $n$ group of data.
\begin{equation}
\label{likelihood}
p(y_1,...,y_n,x_1,...,x_n|\alpha_1,\alpha_2,\tau)=\prod_{i=1}^n\mathcal{N}(\alpha_1x_i^2+\alpha_2x_i+1-\alpha_1-\alpha_2,1/\tau)
\end{equation}
where we place conjugate priors on $\alpha_1, \alpha_2$ and $\tau$:
\begin{eqnarray}
\label{mu1}
\alpha_1&\sim& \mathcal{N}(\mu_1,1/\tau_1)\\
\alpha_2&\sim& \mathcal{N}(\mu_2,1/\tau_2)\\
\tau&\sim& \mathcal{G}(\alpha,\beta)
\end{eqnarray}
\subsubsection{Estimation on $\alpha_1$}
Since only one parameter is updated in each step in GS, updating of all parameters in ZIP model need to follow the procedure introduced in algorithm 2. In order to find out what is the distribution after the updating, the following equation will be considered:
\begin{equation}
\label{posterior}
p(\alpha_1|\alpha_2,\tau,y,x)\propto p(y,x|\alpha_1,\alpha_2,\tau)p(\alpha_1)
\end{equation}
where $p(\alpha_1|\alpha_2,\tau,y,x)$ is the posterior probability, $p(y,x|\alpha_1,\alpha_2,\tau)$ is just the likelihood in \ref{likelihood}, $p(\alpha_1)$ is the prior estimation and it is a normal distribution with mean equals $\mu_1$ , variance equals $1/\tau_1$ in \ref{mu1}. After some calculation, take the log form of \ref{posterior}, the following can be derived:
\begin{equation}
\label{logalpha1}
-\frac{\tau_1}{2}(\alpha_1-\mu_1)^2-\frac{\tau}{2}\sum_{i=1}^n\big(y_i-(\alpha_1x_i^2+\alpha_2x_i+1-\alpha_1-\alpha_2)\big)^2
\end{equation}
The reason for taking log form of this is because it helps convert the multiplications to summations, which can significantly simplify the calculation. Furthermore, the log form helps us update/derive the distribution of posterior more easily and conveniently. For a normal distribution $y$ follows $y\sim \mathcal{N}(\mu,1/\tau)$, the log dependence on $y$ is $-\dfrac{\tau}{2}(y-\mu)^2\propto -\dfrac{\tau}{2}y^2+\tau \mu y$. Therefore, if we can simplify \ref{logalpha1} into that form with variance as $\alpha_1$, the distribution of the posterior will be found. \ref{logalpha1} can be further written as:
\begin{eqnarray}
\nonumber
&-&\frac{\tau_1}{2}\alpha_1^2-2\alpha_1\mu_1+\mu_1^2-\frac{\tau}{2}\sum_{i=1}^n(y_i^2-2y_i\alpha_1x_i^2-2\alpha_2x_iy_i-2y_i+2y_i\alpha_1+2y_i\alpha_2\\
\nonumber
&+&\alpha_1^2x_i^4+\alpha_2x_i^2+1+\alpha_1^2+\alpha_2^2-2\alpha_1-2\alpha_2+2\alpha_1\alpha_2+2\alpha_1\alpha_2x_i^3+2\alpha_1x_i^2\\
\nonumber
&-&2\alpha_1^2x_i^2-2\alpha_1\alpha_2x_i^2+2\alpha_2x_i-2\alpha_1\alpha_2x_i-2\alpha_2^2x_i)
\end{eqnarray}
Because we are interested in finding distribution of $\alpha_1$, all terms not related to $\alpha_1$ can be omitted.
\begin{eqnarray}
\nonumber
&-&\frac{\tau_1}{2}\alpha_1^2-2\alpha_1\mu_1-\frac{\tau}{2}\sum_{i=1}^n(-2y_i\alpha_1x_i^2+2y_i\alpha_1
+\alpha_1^2x_i^4+\alpha_1^2\\
\nonumber
&-&2\alpha_1+2\alpha_1\alpha_2+2\alpha_1\alpha_2x_i^3+2\alpha_1x_i^2
-2\alpha_1^2x_i^2-2\alpha_1\alpha_2x_i^2-2\alpha_1\alpha_2x_i)
\end{eqnarray}
Combining all terms related to $\alpha_1$ and $\alpha_1^2$ yields:
\begin{eqnarray}
\nonumber
-\big(\frac{\tau_1}{2}+\frac{\tau}{2}\sum_{i=1}^n(x_i^2-1)^2\big)\alpha_1^2+\Big(\tau_1\mu_1-\tau\sum_{i=1}^n\big((\alpha_2-1-\alpha_2x_i+y_i)(1-x_i^2)\big)\Big)\alpha_1
\end{eqnarray}
Thus, the posterior follows the following distribution of $\alpha_1$:
\begin{eqnarray}
\nonumber
&\ &\alpha_1|\alpha_2,\tau,\tau_1,\mu_1,x,y,\\
\nonumber
&\ &\sim \mathcal{N}\bigg(\dfrac{\tau_1\mu_1-\tau\sum_{i=1}^n\big((\alpha_2-1-\alpha_2x_i+y_i)(1-x_i^2)\big)}{\tau_1+\tau\sum_{i=1}^n(x_i^2-1)^2},\dfrac{1}{\tau_1+\tau\sum_{i=1}^n(x_i^2-1)^2}\bigg)
\end{eqnarray}
\subsubsection{Estimation on $\alpha_2$}
Same procedure can be applied to derive conditional distribution of $\alpha_2$:
\begin{eqnarray}
\nonumber
&\ &\alpha_2|\alpha_1,\tau,\tau_2,\mu_2,x,y,\\
\nonumber
&\ &\sim \mathcal{N}\bigg(\dfrac{\tau_2\mu_2-\tau\sum_{i=1}^n\big((\alpha_1-1-\alpha_1x_i^2+y_i)(1-x_i)\big)}{\tau_2+\tau\sum_{i=1}^n(x_i-1)^2},\dfrac{1}{\tau_2+\tau\sum_{i=1}^n(x_i-1)^2}\bigg)
\end{eqnarray}
\subsubsection{Estimation on $\tau$}
\begin{figure}[!t]
\centering
    \includegraphics[width=\linewidth,]{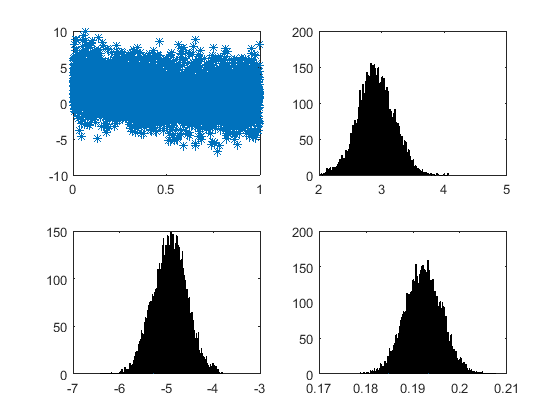}
\caption{Gibbs Sampling example .}
\label{ZIP1}
\end{figure}

The Gamma density is given by $p(x;\alpha,\beta)\propto \beta^{\alpha}x^{\alpha-1}e^{-\beta x}$, its log form is proportional  to $(\alpha-1)\log x-\beta x$. The posterior can be written as $p(\tau|\alpha_1,\alpha_2,y,x)\propto p(y,x|\alpha_1,\alpha_2,\tau)p(\tau)$. Using same procedure, we can get the log expression as:
\begin{equation}
\nonumber
\frac{n}{2}\log \tau-\frac{\tau}{2}\sum_{i=1}^n\big(y_i-\alpha_1x_i^2-\alpha_2x_i-(1-\alpha_1-\alpha_2)\big)+(\alpha-1)\log \tau-\beta \tau
\end{equation}
Therefore, distribution of $\tau$ follows Gamma distribution:
\begin{eqnarray}
\nonumber
&\ &\tau|\alpha_1,\alpha_2,\alpha,\beta,x,y,\\
\nonumber
&\ &\sim \mathcal{G}\bigg(\alpha+\dfrac{n}{2},\beta+\dfrac{\sum_{i=1}^n\big(y_i-\alpha_1x_i^2-\alpha_2x_i-(1-\alpha_1-\alpha_2)\big)^2}{2}\bigg)
\end{eqnarray}
Simulation verified this two parameters estimation. The result is shown in  Fig. \ref{ZIP1}. The true distribution is $\alpha_1=3$ (upper right corner), $\alpha_2=-5$ (lower left corner), the precision, which equals one over variance is 0.2 (shown in the lower right corner). The theoretical results and the Gibbs sampling results are compared in Table. \ref{table1}.\par
\begin{table}[!t]
\caption{Comparison of Gibbs and theoretical results}
\centering
\label{table1}
\begin{tabular}{ c | c | c }
  \hline
  \centering			
  Parameter & Theoretical Result & Gibbs Sampling Result\\
  \centering
  $\alpha_1$& 3 & 2.98 \\
 $\alpha_2$ & -5 &-4.99 \\
$\alpha_3$ & 3 & 3.02 \\
$\tau$ & 0.2 & 0.192 \\
  \hline  
\end{tabular}
\end{table}

Same procedures were followed to derive three parameters ($\alpha_1, \alpha_2, \alpha_3$) estimation. Updated distributions of the three parameters are:
\begin{eqnarray}
\nonumber
&\ &\alpha_1|\alpha_2,\alpha_3, \tau,\tau_1,\mu_1,x,y,\\
\nonumber
&\ &\sim \mathcal{N}\bigg(\dfrac{\tau_1\mu_1+\tau\sum_{i=1}^n\big((y_i-\alpha_3)x_i^2-\alpha_2x_i^3\big)}{\tau_1+\tau\sum_{i=1}^n(x_i)^4},\dfrac{1}{\tau_1+\tau\sum_{i=1}^n(x_i)^4}\bigg)
\end{eqnarray}
\begin{eqnarray}
\nonumber
&\ &\alpha_2|\alpha_1,\alpha_3, \tau,\tau_2,\mu_2,x,y,\\
\nonumber
&\ &\sim \mathcal{N}\bigg(\dfrac{\tau_2\mu_2+\tau\sum_{i=1}^n\big((y_i-\alpha_3)x_i-\alpha_1x_i^3\big)}{\tau_2+\tau\sum_{i=1}^n(x_i)^2},\dfrac{1}{\tau_2+\tau\sum_{i=1}^n(x_i)^2}\bigg)
\end{eqnarray}
\begin{eqnarray}
\nonumber
&\ &\alpha_3|\alpha_1,\alpha_2, \tau,\tau_3,\mu_3,x,y,\\
\nonumber
&\ &\sim \mathcal{N}\bigg(\dfrac{\tau_3\mu_3+\tau\sum_{i=1}^n\big(y_i-\alpha_1x_i^2-\alpha_2x_i\big)}{\tau_3+\tau\sum_{i=1}^n(x_i)^0},\dfrac{1}{\tau_3+\tau\sum_{i=1}^n(x_i)^0}\bigg)
\end{eqnarray}
\begin{eqnarray}
\nonumber
&\ &\tau|\alpha_1,\alpha_2,,\alpha_3\, \alpha,\beta,x,y,\\
\nonumber
&\ &\sim \mathcal{G}\bigg(\alpha+\dfrac{n}{2},\beta+\dfrac{\sum_{i=1}^n\big(y_i-\alpha_1x_i^2-\alpha_2x_i-\alpha_3\big)^2}{2}\bigg)
\end{eqnarray}

\begin{figure}[!t]
\centering
    \includegraphics[width=1\linewidth,]{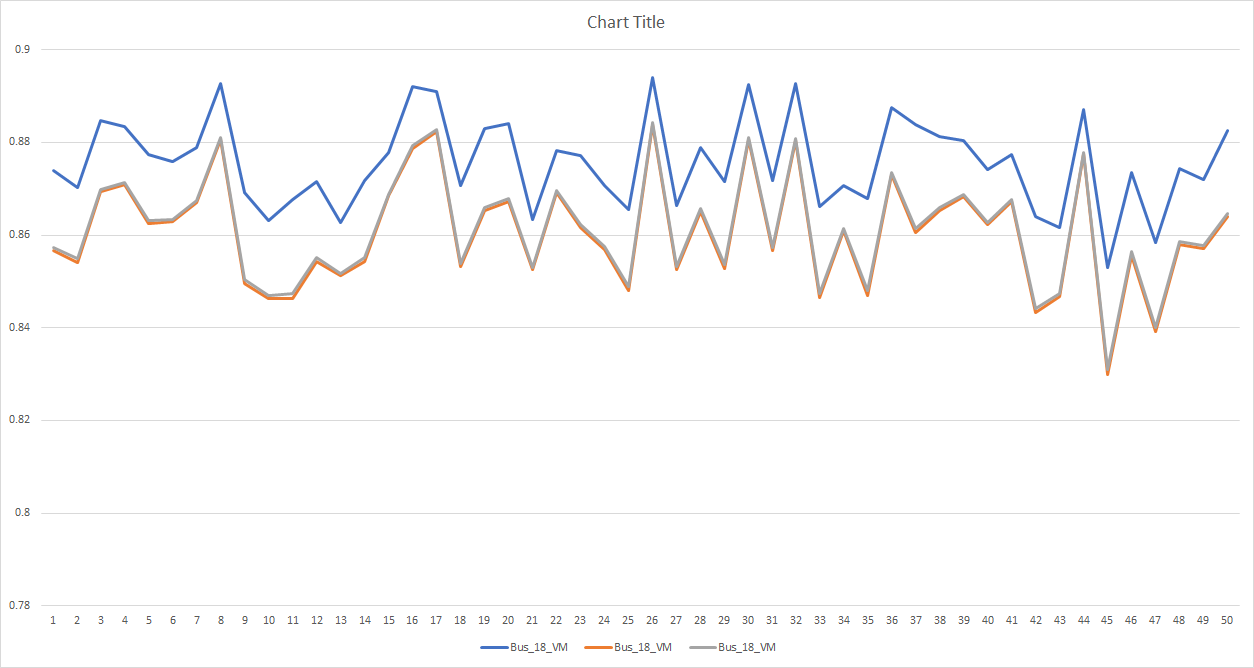}
\caption{Comparison 1 using random normal distributed load at all 33 buses, plotting the voltage at bus 18 using both the real parameter and estimated parameter. Blue line at the top is the voltage at bus 18 without adding the ZIP model.}
\label{comp1}
\end{figure}
\begin{figure}[!t]
\centering
    \includegraphics[width=1\linewidth,]{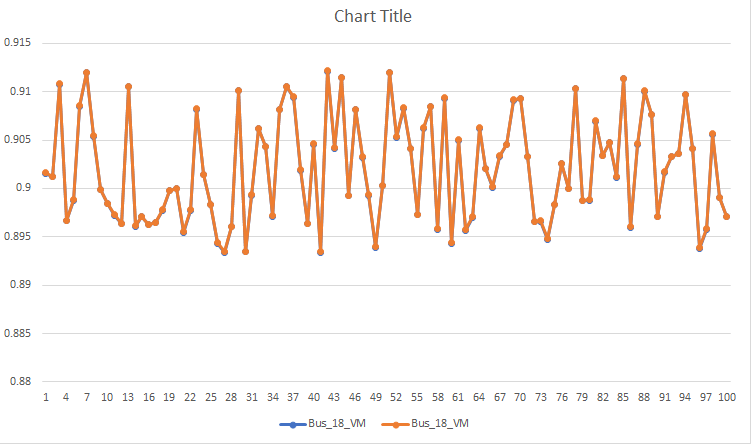}
\caption{Comparison 1 using random normal distributed load in 33 bus test feeder, plotting the voltage at bus 18 using both the real parameter and estimated parameter.}
\label{comp2}
\end{figure}
\begin{figure}[!t]
\centering
    \includegraphics[width=1\linewidth,]{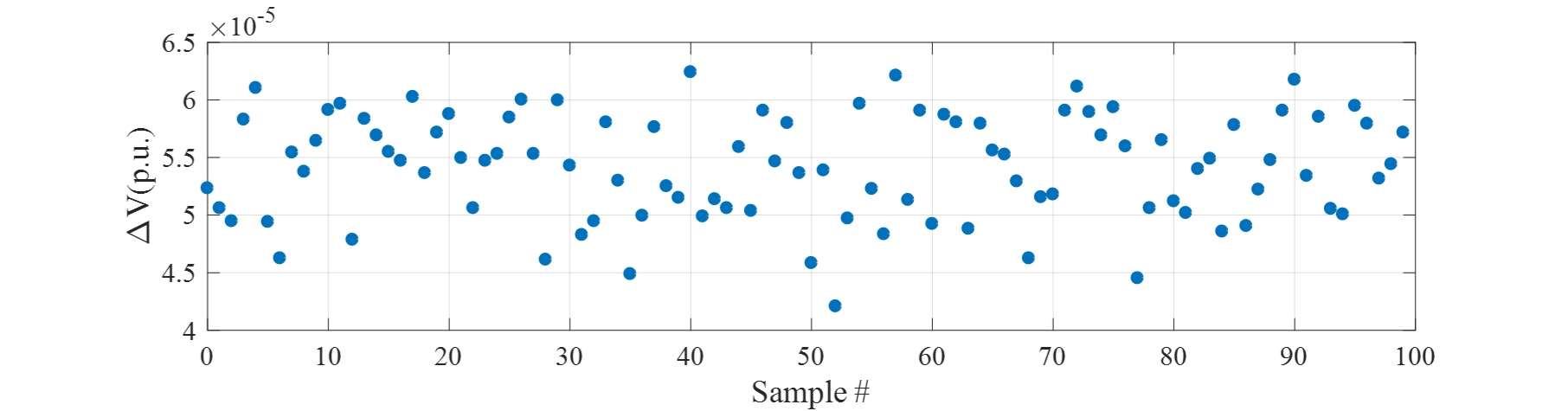}
\caption{Comparison 1 of $\Delta V$ between connecting real ZIP model and the estimated ZIP model at bus 18.}
\label{comp3}
\end{figure}
\begin{figure}[!t]
\centering
    \includegraphics[width=1\linewidth,]{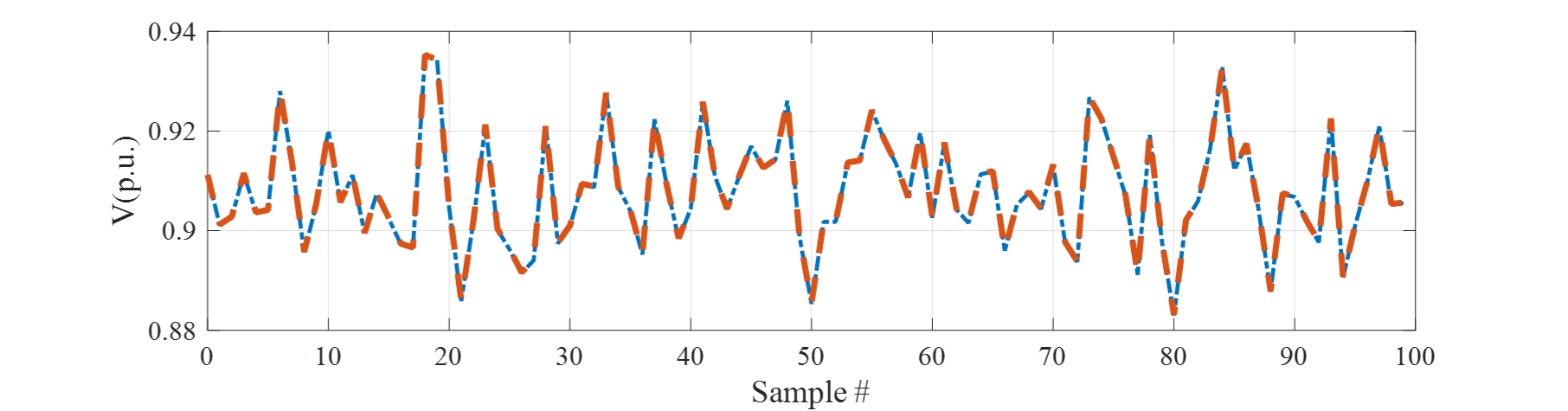}
\caption{Comparison 2 in the 33 bus system, plotting the voltage at bus 18 using both the real parameter and estimated parameter.}
\label{comp4}
\end{figure}
\begin{figure}[!t]
\centering
    \includegraphics[width=1\linewidth,]{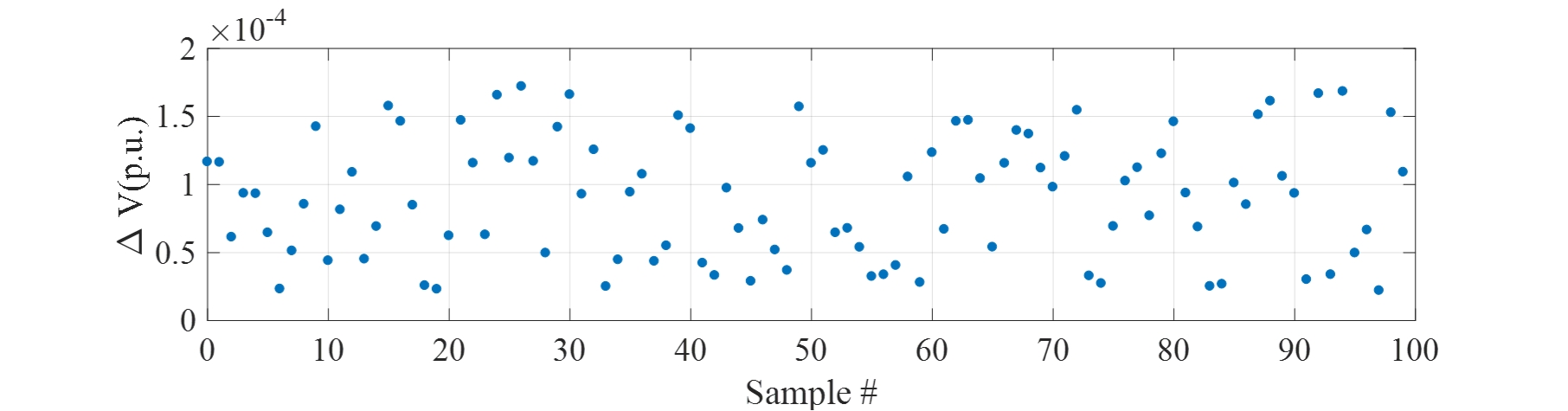}
\caption{Comparison 2 of $\Delta V$ between connecting real ZIP model and the estimated ZIP model at bus 18.}
\label{comp5}
\end{figure}
\begin{figure}[!t]
\centering
    \includegraphics[width=1\linewidth,]{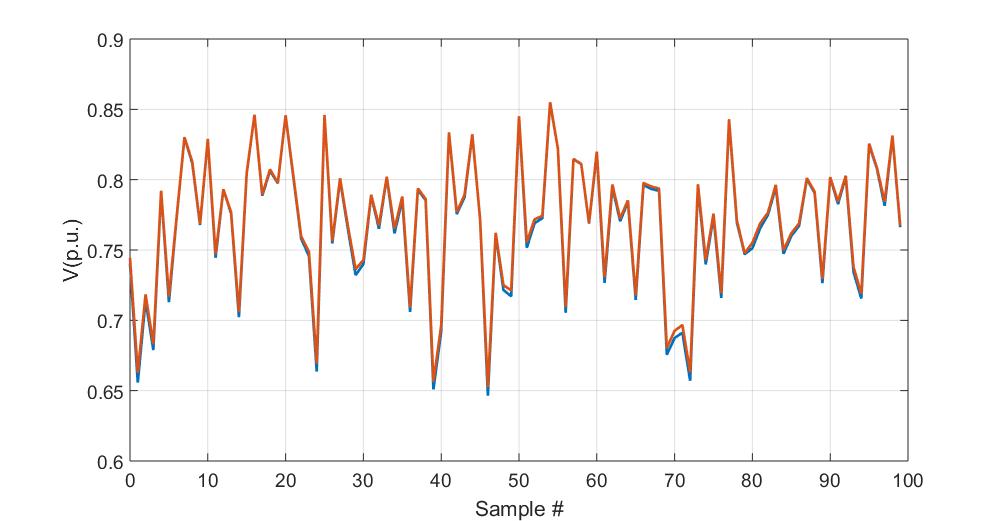}
\caption{Comparison 3 in the 33 bus system, plotting the voltage at bus 18 using both the real parameter and estimated parameter.}
\label{comp6}
\end{figure}
\begin{figure}[!t]
\centering
    \includegraphics[width=1\linewidth,]{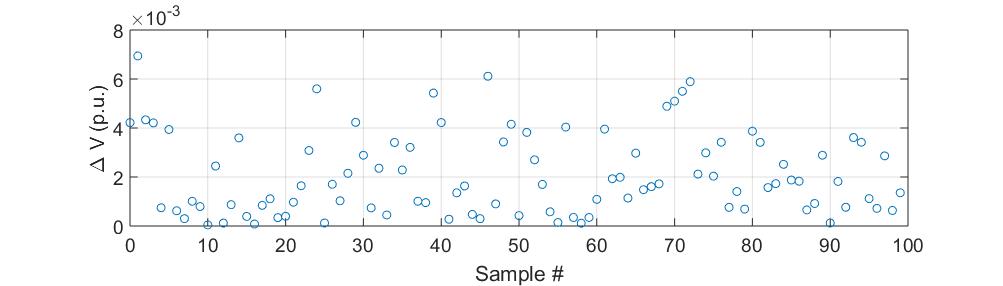}
\caption{Comparison 3 of $\Delta V$ between connecting real ZIP model and the estimated ZIP model at bus 18.}
\label{comp7}
\end{figure}
\begin{figure}[!t]
\centering
    \includegraphics[width=1\linewidth,]{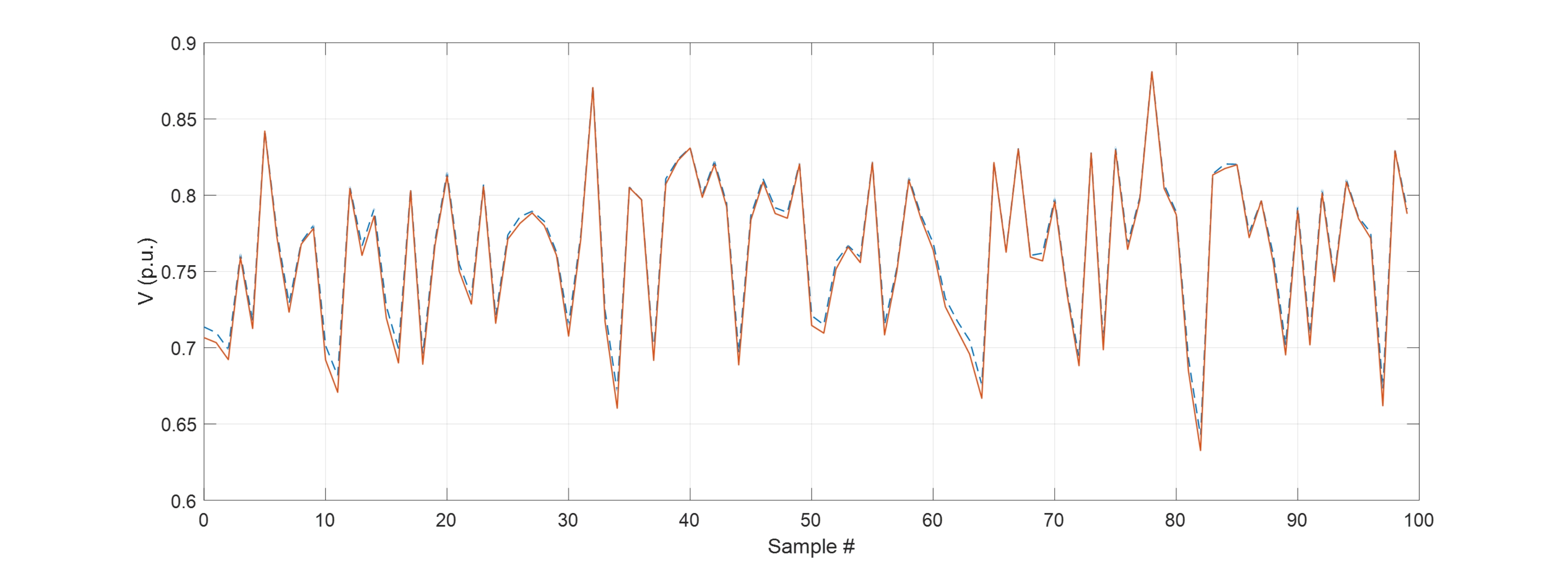}
\caption{Comparison in 33 bus system, plotting the voltage at bus 18 using both the real parameter and estimated parameter.}
\label{comp8}
\end{figure}
\begin{figure}[!t]
\centering
    \includegraphics[width=1\linewidth,]{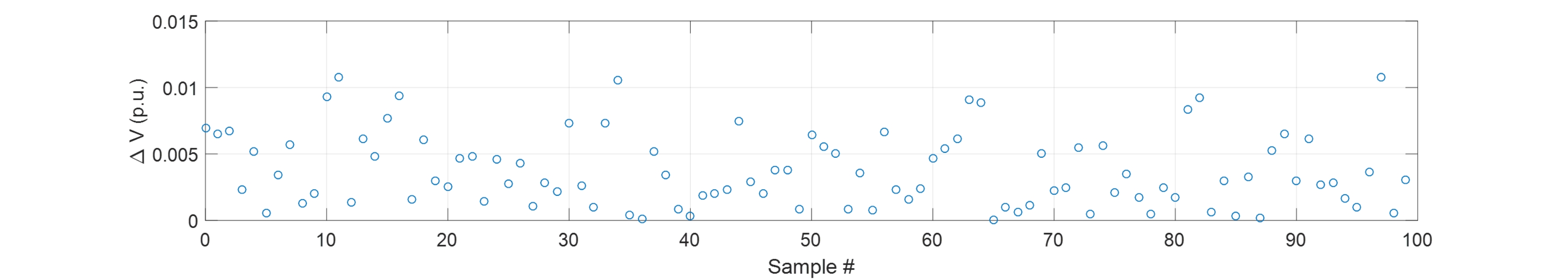}
\caption{Comparison of $\Delta V$ between connecting real ZIP model and the estimated ZIP model at bus 18 using an arbitrary ZIP model and the real ZIP model.}
\label{comp9}
\end{figure}
\subsubsection{Gibbs Sampling Experiment for ZIP Model}
In order to verify the proposed parameter estimation method, a 33 bus distribution test feeder connected with a fixed ZIP model was used to generate the simulation data. Three experiments were made: a random coefficient follows $\mathcal{N}(2.0.8)$ , $\mathcal{U}[0.3,1.8]$, $\mathcal{U}{[0.01,4.5]}$ is multiplied at each load in the 33 bus distribution system in each experiment, respectively. Differences between the three experiments are the changing ranges of the voltage. In the first experiment and the second experiment, the changing ranges are 0.02 p.u. and 0.06 p.u., while in the last one, the range was extended to 0.2 p.u.  \par
 Bus 18 was selected as the measurement bus. A parameter fixed ZIP model was connected to this bus. In each experiment, 1000 independent  power flow calculations were made by using PSAT. Coefficients draw from the aforementioned distributions were multiplied at each load in all 33 buses to simulate the variance of the load. Voltage ($V$) at bus 18 was measured in every iteration, as well as the active power $P$. The initial voltage $V_0$ and active power $P_0$ at bus 18 was also measure because in each iteration, different multiplier will change the voltage and active power at the corresponding bus at each node. \par 
\begin{figure}[!t]
\centering
    \includegraphics[width=1\linewidth,]{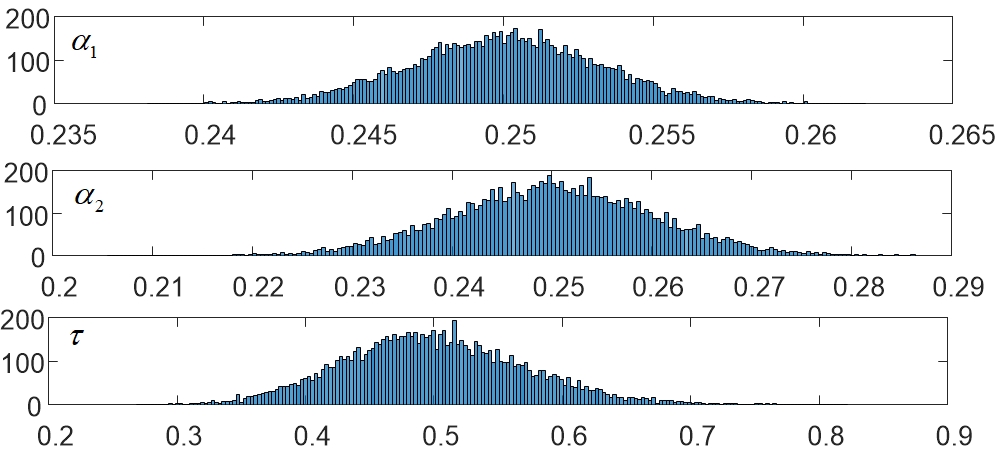}
\caption{The estimated parameters of the ZIP model.}
\label{para1}
\end{figure}
\subsubsection{Comparison Results}
\begin{figure}[!t]
\centering
    \includegraphics[width=1\linewidth,]{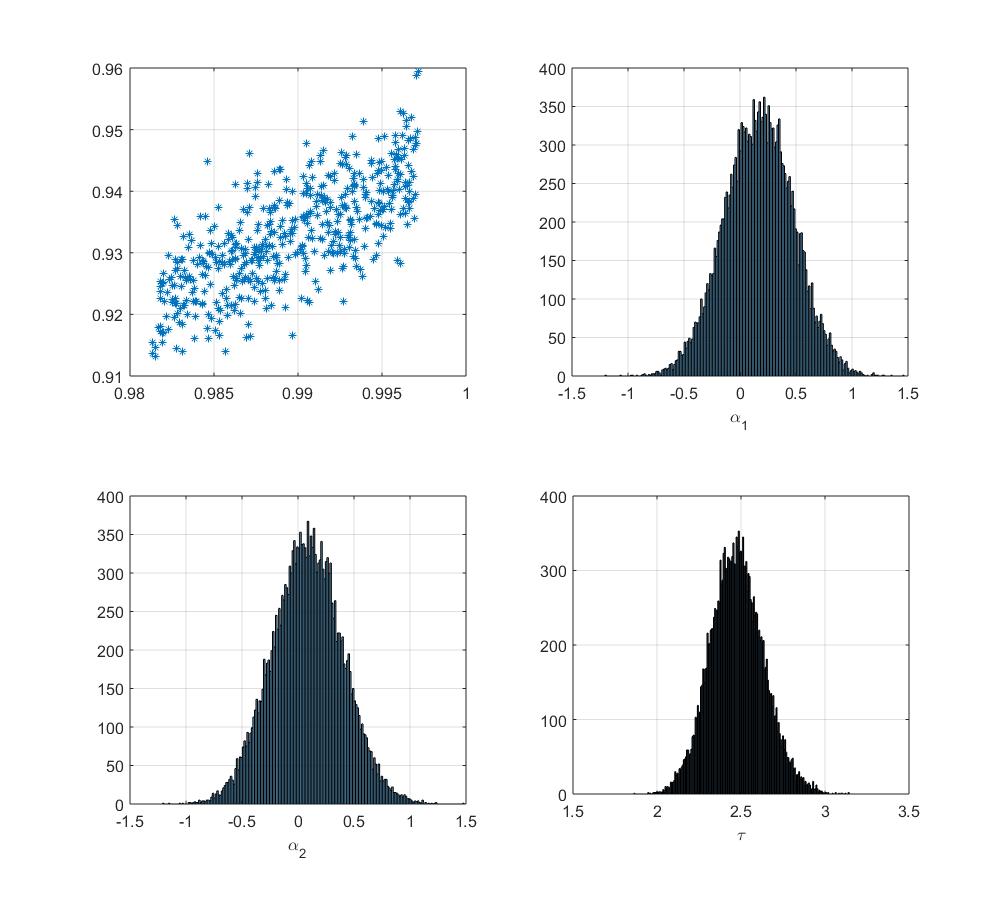}
\caption{The estimated parameters of the ZIP model.}
\label{para2}
\end{figure}

The estimations of using GS for two parameters and three parameters are shown in Figs. \ref{para2} and \ref{para1}, respectively. The real parameters of the ZIP model is $\alpha_1=0.25, \alpha_2=0.25, \alpha_3=0.5$. The mean value for $\alpha_1$ and $\alpha_2$ are 0.21 and 0.27, respectively, in Fig. \ref{para2}. Consequently, the last parameter $\alpha_3$ equals 1-0.21-0.17=0.62.\par
For the three parameters estimation, the sampled parameters is shown in Fig. \ref{para1}. The mean of the estimated three parameters are $\alpha_1=0.31$, $\alpha_2=0.32$,$\alpha_3=0.31$. However, recall that $\alpha_1+\alpha_2+\alpha_3=1$, therefore, the three parameters are chosen as: $\alpha_1=0.31, \alpha_2=0.32, \alpha_3=0.37$. The estimated ZIP model is then compared with the real ZIP model by run the simulation with same load parameters. 
Three comparisons of the proposed GS parameter estimations are shown in Figs. \ref{comp1} to \ref{comp7} since the estimated parameters are not as close as the real values. These comparisons are going to show the impact of the parameters of the load in terms of different voltage changes. \par
As they were shown in Figs. \ref{comp2}, \ref{comp4} and \ref{comp6}, the estimated ZIP model has almost the same voltage response compared with the real ZIP model. The two lines in Figs. \ref{comp2}, \ref{comp4} and \ref{comp6} coincide with each other and the differences $\Delta V$ in three comparisons are shown in Figs.\ref{comp3}, \ref{comp5}, \ref{comp7}, respectively. It is clearly shown that the differences of the voltage at bus 18 is relatively small. The magnitude of $\Delta V$ only reaches $10^{-5}$ when the load change is small. With the increasing changing of the load, in Fig. \ref{comp6}, when the voltage floating within [0.65 p.u., 0.88 p.u.], the voltage response for the ZIP model only differed by $10^{-3}$.  According to the three comparisons, it is indicated the estimated parameters are suitable for usage. \par
In Figs. \ref{comp8} and \ref{comp9}, the voltage response of an arbitrary selected ZIP coefficients were compared with the real 
ZIP model. It is shown that $\Delta V$ increased to $10_{-2}$, which is 10 times larger than that in Fig. \ref{comp7}. So there is necessarity to identify the coefficients in a certain range.\par
%

\section{Induction Motor (IM)}
\subsection{Model of IM}

\begin{figure}[!t]
\centering
    \includegraphics[width=0.75\linewidth,height=0.87\textheight]{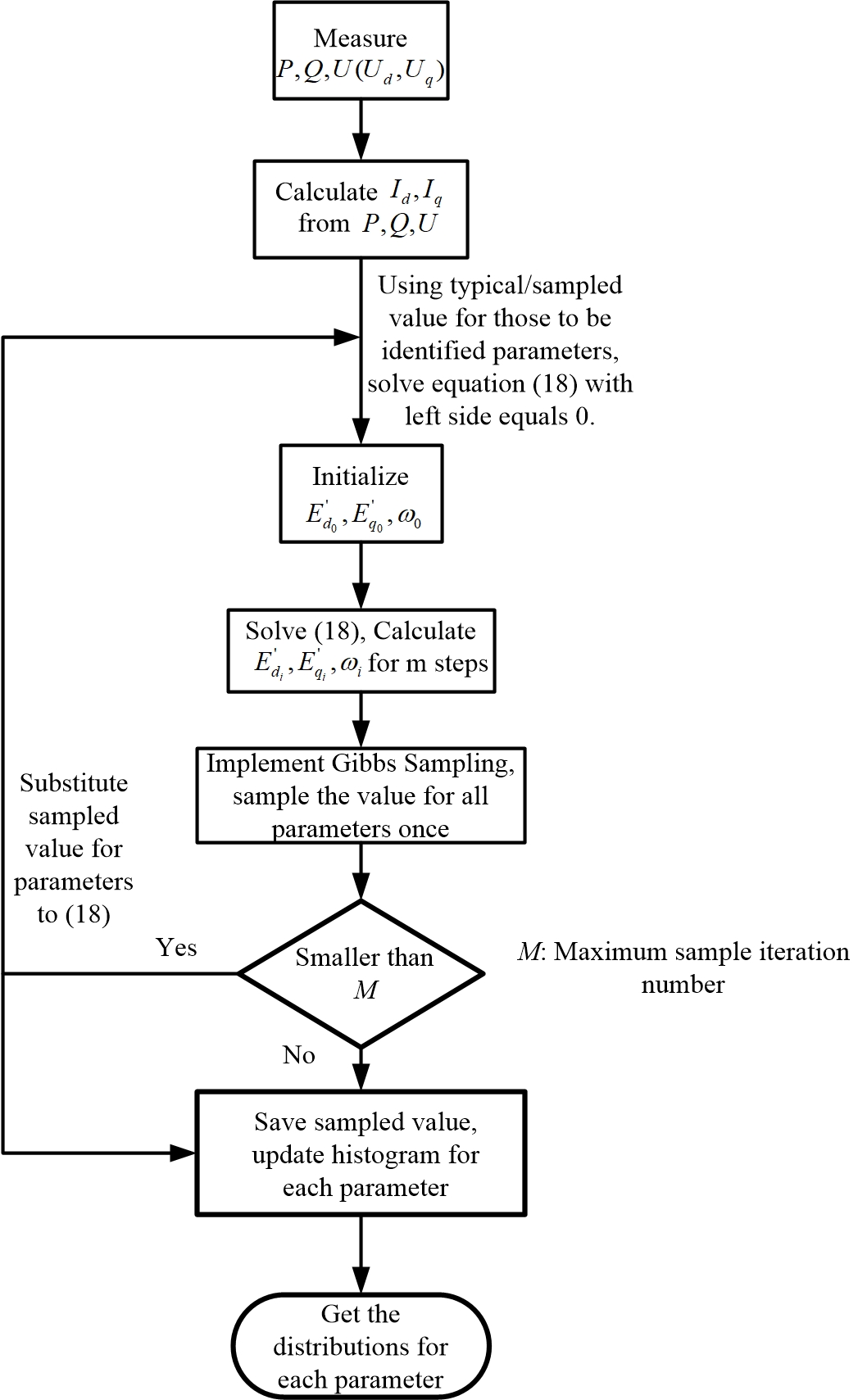}
\caption{Flow chart for IM model using GS method.}
\label{imgibbs}
\end{figure}
The dynamic part of a composite load is contributed by induction motor (IM) in this study. Detailed model of a IM has been introduced in many lectures and papers. In this study, we selected the model used in \cite{gibbsmjin}. The dynamic model of IM can be written as follows:
\begin{equation}
\label{im1}
\begin{cases}
\dfrac{dE'_d}{dt}=-\dfrac{1}{T'}[E'_d+(X-X')I_q]-(\omega-1)E'_q\\
\dfrac{dE'_q}{dt}=-\dfrac{1}{T'}[E'_q-(X-X')I_d]+(\omega-1)E'_d\\
\dfrac{d\omega}{dt}=-\dfrac{1}{2H}[(A^2\omega+B\omega+C)T_0-(E'_dI_d+E'_qI_q)]
\end{cases}
\end{equation}
\begin{equation}
\label{im2}
\begin{cases}
I_d=\dfrac{1}{R_s^2+X'^2}[R_s(U_d-E'_d)+X'(U_q-E'_q)]\\
I_q=\dfrac{1}{R_s^2+X'^2}[R_s(U_q-E'_q)-X'(U_d-E'_d)]
\end{cases}
\end{equation}
where 
\begin{eqnarray}
\nonumber
T'&=&\dfrac{X_r+X_m}{R_r}\\
\nonumber
X&=&X_s+X_m\\
\nonumber
X'&=&X_s+\dfrac{X_mX_r}{X_m+X_r}\\
\nonumber
A&+&B\ +\ C=1
\end{eqnarray}
where $R_s$ is the stator winding resistance of motor, $X_s$ is  stator leakage reactance of motor, $X_m$ is  magnetizing reactance of motor,$R_r$ is  rotor resistance, $X_r$  is  rotor leakage reactance, $H$ is  rotor inertia constant, $\omega $ is  rotor speed, $I_d, I_q$ are stator current in $d$-axis and $q$-axis, $U_d, U_q$ are  bus voltage in $d$-axis and $q$-axis, $E_d, E_q$ are rotor flux linkages in $d$-axis and $q$-axis.\par
\begin{figure}[!h]
\centering
    \includegraphics[width=\linewidth,height=0.65\textheight]{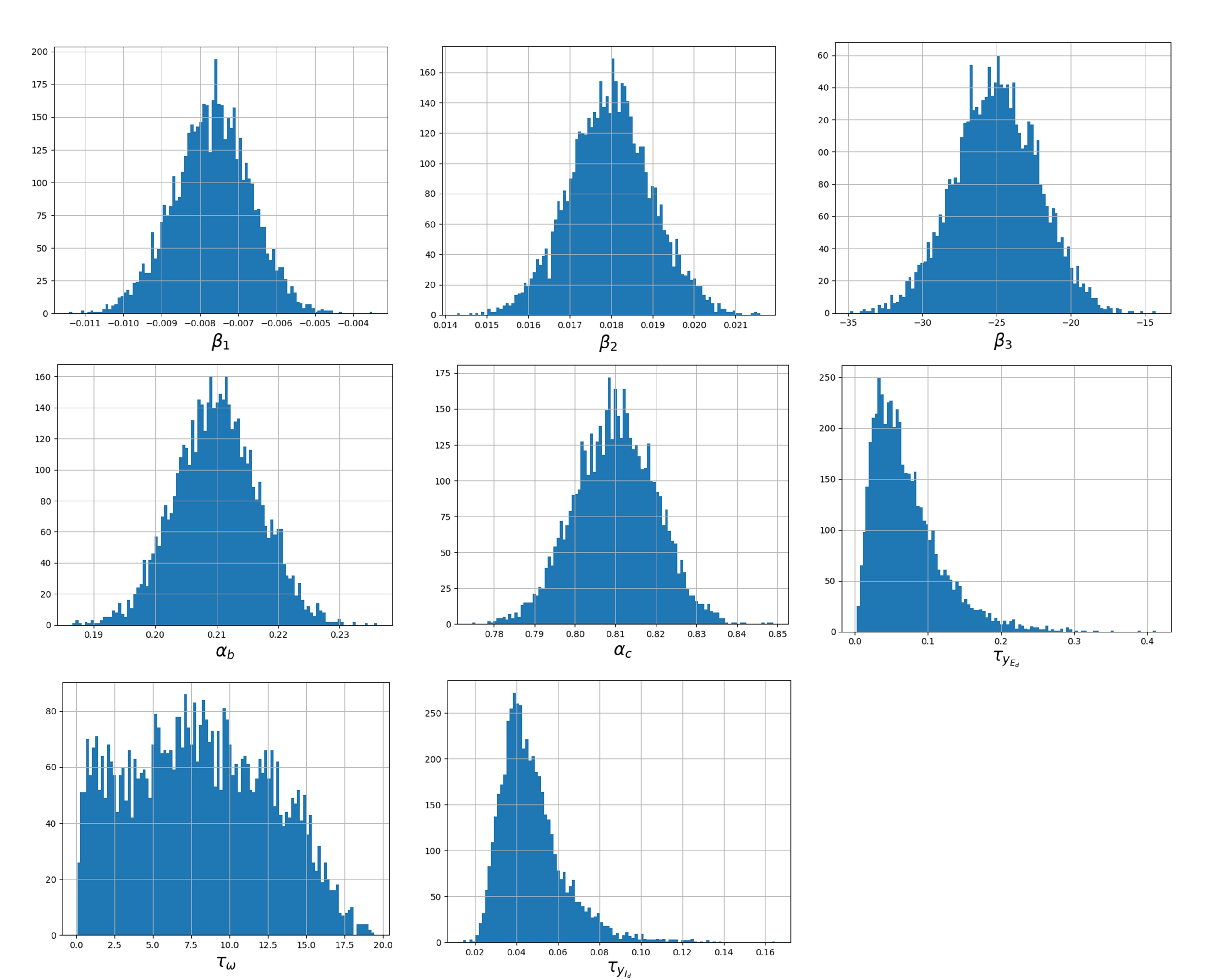}
\caption{The estimated parameters of the ZIP model.}
\label{imgibbs}
\end{figure}
Differences of the coefficients in first two equations in \ref{im1} are only the signs. Therefore, only one of those two equations will be used in the estimation. Mechanical torque is assumed to be proportional to the square of the rotation speed of the IM, thus the value of A is 1 while B and C are all 0. Similarly, in \ref{im2}, both of the two equations can be used to estimate $R$ and $X'$, and we will propose Gibbs sampling by using the equation of $I_d$.
As introduced in previous section, Gibbs sampling requires update all the parameters iteratively in each step. Define $y_{Ed} \triangleq \dfrac{dE'_q}{dt}$, $y_{Eq} \triangleq \dfrac{dE'_d}{dt}$ $y_{\omega}\triangleq\dfrac{d\omega}{dt}$, $ y_{Id}\triangleq I_d$, $ y_{Iq}\triangleq I_q$, $ \beta_1\triangleq-\dfrac{1}{T'}$, $\beta_2\triangleq-\dfrac{X-X'}{T'}$, $\beta_3\triangleq-\dfrac{1}{2H}$, $\alpha_b \triangleq \dfrac{R_s}{R_s^2+X'^2}$, $\alpha_c \triangleq\dfrac{X'}{R_s^2+X'^2}$. the equations will be used in the GS can be written as:

\begin{eqnarray}
\label{equa1}
y_{Ed}[t]=\beta_1E'_d[t]+\beta_2I_q[t]-(\omega-1)E'_q[t]+\varepsilon_{E}[t]\\
\label{equa2}
y_{Eq}[t]=\beta_1E'_q[t]-\beta_2I_d[t]+(\omega-1)E'_d[t]+\varepsilon_{E}[t]\\
\label{equa3}
y_{\omega}[t]=\beta_3(\omega^2T_0-E'_d[t]I_d[t]-E'_q[t]I_q[t])+\varepsilon_{\omega}[t]\\
\label{equa4}
y_{Id}[t]=\alpha_b(U_d[t]-E'_d[t])+\alpha_c(U_q[t]-E'_q[t])+\varepsilon_{I}[t]\\
\label{equa5}
y_{Iq}[t]=\alpha_b(U_q[t]-E'_q[t])-\alpha_c(U_d[t]-E'_d[t])+\varepsilon_{I}[t]
\end{eqnarray}
where ${\varepsilon_{E}[t]\sim \mathcal{N}(0,1/\tau_{E})}$, ${\varepsilon_{\omega}[t]\sim \mathcal{N}(0,1/\tau_{\omega})}$, and ${\varepsilon_{I}[t]\sim \mathcal{N}(0,1/\tau_{I})}$.
\begin{eqnarray}
\nonumber
\beta_1&\sim& \mathcal{N}(\mu_{\beta1}^{(0)},1/\tau_{\beta1}^{(0)}), \beta_2\sim \mathcal{N}(\mu_{\beta2}^{(0)},1/\tau_{\beta2}^{(0)})\\
\nonumber
\beta_3&\sim& \mathcal{N}(\mu_{\beta3}^{(0)},1/\tau_{\beta3}^{(0)}), \alpha_b\sim \mathcal{N}(\mu_{\alpha b}^{(0)},1/\tau_{\alpha b}^{(0)})\\
\nonumber
\alpha_c&\sim& \mathcal{N}(\mu_{\alpha c}^{(0)},1/\tau_{\alpha c}^{(0)}), \tau_{{E}}\sim \mathcal{G}(\alpha_{{E}}^{(0)},\beta_{{E}}^{(0)})\\
\nonumber
\tau_{{\omega}}&\sim& \mathcal{G}(\alpha_{{\omega}}^{(0)},\beta_{{\omega}}^{(0)}), \ \tau_{{I}}\sim \mathcal{G}(\alpha_{{I}}^{(0)},\beta_{{I}}^{(0)}).
\end{eqnarray}
Similar as in previous section, the updated conditional distributions of each parameters given by i.i.d experiments are derived by using the same approach. Representations of these five parameters can be found as follows:

Define
\[
\begin{array}{l}
\mu_{\beta_1}^{(1)}=\\
~~~~\dfrac{\tau_{\beta_1}^{(0)}\mu_{\beta_1}^{(0)}+\tau_{E}^{(0)}\sum_{t=1}^n\big(E'_{d}[t]y_{E_d}[t]+E'_{q}[t]y_{E_q}[t]-\beta_2^{(0)}(E'_{d}[t]I_{q}[t]-E'_{q}[t]I_{d}[t])\big)}{\tau_1^{(0)}+\tau_{E_d}^{(0)}\sum_{t=1}^n(E_{d}'[t]^{2}+E_{q}'[t]^{2})},\\
\tau_{\beta_1}^{(1)}=\dfrac{1}{\tau_{\beta_1}^{(0)}+\tau_{E}^{(0)}\sum_{t=1}^n(E_{d}'[t]^{2}+E_{q}'[t]^{2})}
\end{array}
\]
the posterior  $\beta_1^{(1)}|\  \beta_2^{(0)}, y_{E_d}, y_{E_q}, E'_d, E'_q, I_q, I_d, \omega, \tau_{E}^{(0)},\tau_{\beta_1}^{(0)},\mu_{\beta _1}^{(0)}\sim \mathcal{N}(\mu_{\beta_1}^{(1)},\tau_{\beta_1}^{(1)})$ can be obtained.\par
Define
\[\small
\begin{array}{l}
\mu_{\beta_2}^{(1)}=\\
~~~~\dfrac{\tau_{\beta_2}^{(0)}\mu_{\beta_2}^{(0)}+\tau_{E}^{(0)}\sum_{t=1}^n\big(y_{E_d}[t]I_{q}[t]-y_{E_q}[t]I_{d}[t]+\beta_1^{(1)}(E'_{d}[t]I_{q}[t]-E'_{q}[t]I_{d}[t])+(\omega[t]-1)(E'_{d}[t]I_{d}[t]+E_{q}[t]I_{q}[t])\big)}{\tau_{\beta_ 1}^{(0)}+\tau_{E}^{(0)}\sum_{t=1}^n(I_{d}[t]^{'2}+I_{q}[t]^{'2})},\\
\tau_{\beta_2}^{(1)}=\dfrac{1}{\tau_{\beta_2}^{(0)}+\tau_{E}^{(0)}\sum_{i=1}^n(I_{d}'[t]^{2}+I_{q}'[t]^{2})}
\end{array}
\]
\normalsize
the posterior  $\beta_2^{(1)}|\  \beta_1^{(1)}, y_{E_d},y_{E_q}, E'_d, E'_q, I_q, I_d, \omega, \tau_{E}^{(0)},\tau_{\beta_2}^{(0)},\mu_{\beta_2}^{(0)}\sim \mathcal{N}(\mu_{\beta_2}^{(1)},\tau_{\beta_2}^{(1)})$ can be obtained.\par

Define
\[
\begin{array}{l}
\mu_{\beta_3}^{(1)}=\dfrac{\tau_{\beta_3}^{(0)}\mu_{\beta_3}^{(0)}-\tau_{\omega}^{(0)}\sum_{t=1}^n\big(y_{\omega}[t](-\omega[t]^2T_0+E'_{d}[t]I_{d}[t]+E_{q}[t]I_{q}[t])\big)}{\tau_{\beta 3}^{(0)}+\tau_{\omega}^{(0)}\sum_{t=1}^n\big(\omega[t]^4T_0^2+(E'_{d}[t]I_{d}[t]+E_{q}[t]I_{q}[t])^2-2\omega[t]^2T_0(E'_{d}[t]I_{d}[t]+E_{q}[t]I_{q}[t])\big)},\\
\tau_{\beta_3}^{(1)}=\dfrac{1}{\tau_{\beta_ 3}^{(0)}+\tau_{\omega}^{(0)}\sum_{t=1}^n\big(\omega[t]^4T_0^2+(E'_{d}[t]I_{d}[t]+E_{q}[t]I_{q}[t])^2-2\omega[t]^2T_0(E'_{d}[t]I_{d}[t]+E_{q}[t]I_{q}[t])\big)}
\end{array}
\]
the posterior  $\beta_3^{(1)}| y_{\omega}, E'_d, E'_q, I_q, I_d, \omega, \tau_{E}^{(0)},\tau_{\beta_3}^{(0)},\mu_{\beta_3}^{(0)}\sim \mathcal{N}(\mu_{\beta_3}^{(1)},\tau_{\beta_3}^{(1)})$ can be obtained.\par

Define
\[
\begin{array}{l}
\mu_{\alpha_b}^{(1)}=\dfrac{\tau_{\alpha_b}^{(0)}\mu_{\alpha_b}^{(0)}+\tau_{I}^{(0)}\sum_{i=1}^n\big(I_{di}(U_{di}-E'_{di})+I_{qi}(U_{qi}-E'_{qi})\big)}{\tau_{\alpha_b}^{(0)}+\tau_{I}^{(0)}\sum_{i=1}^n\big((U_{di}-E_{di})^2+(U_{qi}-E_{qi})^2\big)},\\
\tau_{\alpha_b}^{(1)}=\dfrac{1}{\tau_{\alpha_b}^{(0)}+\tau_{I}^{(0)}\sum_{i=1}^n\big((U_{di}-E_{di})^2+(U_{qi}-E_{qi})^2\big)}
\end{array}
\]
the posterior  $\alpha_b^{(1)}|\alpha_c^{(0)}, y_{Id},y_{Iq}, E'_d, E'_q, I_d, \omega, \tau_{I}^{(0)},\tau_{\alpha_b}^{(0)},\mu_{\alpha_b}^{(0)}\sim \mathcal{N}(\mu_{\alpha_b}^{(1)},\tau_{\alpha_b}^{(1)})$ can be obtained.\par

Define
\[
\begin{array}{l}
\mu_{\alpha_c}^{(1)}=\dfrac{\tau_{\alpha_c}^{(0)}\mu_{\alpha_c}^{(0)}+\tau_{I}^{(0)}\sum_{i=1}^n\big(I_{di}(U_{di}-E'_{di})-I_{qi}(U_{qi}-E'_{qi})\big)}{\tau_{\alpha_c}^{(0)}+\tau_{I}^{(0)}\sum_{i=1}^n\big((U_{di}-E_{di})^2+(U_{qi}-E_{qi})^2\big)},\\
\tau_{\alpha_c}^{(1)}=\dfrac{1}{\tau_{\alpha_c}^{(0)}+\tau_{I}^{(0)}\sum_{i=1}^n\big((U_{di}-E_{di})^2+(U_{qi}-E_{qi})^2\big)}
\end{array}
\]
the posterior  $\alpha_c^{(1)}|\alpha_b^{(1)}, y_{Id},y_{Iq}, E'_d, E'_q, I_d, \omega, \tau_{I}^{(0)},\tau_{\alpha_c}^{(0)},\mu_{\alpha_c}^{(0)}\sim \mathcal{N}(\mu_{\alpha_c}^{(1)},\tau_{\alpha_c}^{(1)})$ can be obtained.\par
Other parameters can be obtained as follows:
\begin{eqnarray}
\nonumber
&\ & \tau_{E}^{(1)}|y_{Ed},y_{Eq}, I_q, E'_d, \omega, E'_q, I_d,\beta_1^{(1)},\beta_2^{(1)}\\
\nonumber
&\ &\sim \mathcal{G}\bigg(\alpha_{E}^{(0)}+\dfrac{n}{2}, \\
\nonumber
&\ &\beta_{E}^{(0)}+\dfrac{\sum_{i=1}^{n_1}\big((y_{Ed}-\beta_1^{(1)}E'_{qi}-\beta_2^{(1)}I_{di}+(\omega_i-1)E'_{qi})^2+(y_{Eq}-\beta_1^{(1)}E'_{di}+\beta_2^{(1)}I_{qi}-(\omega_i-1)E'_{di})^2\big)}{2}\bigg)\\
\nonumber
&\ & \tau_{\omega}^{(1)}|y_{\omega}, I_q,E'_d,T_0,E'_q,U_q,U_d,I_d, \omega, \beta_3^{(1)}\\
\nonumber
&\ &\sim \mathcal{G}\bigg(\alpha_{\omega}^{(0)}+\dfrac{n}{2}, \beta_{\omega}^{(0)}+\dfrac{\sum_{i=1}^{n_2}\big(y_{\omega i}-\beta_3^{(1)}\big(\omega_i^2T_0-(E'_{di}I_{di}+E'_{qi}I_{qi})\big)^2}{2}\bigg)\\
\nonumber
&\ & \tau_{I}^{(1)}|y_{Id},y_{Iq}, I_q,E'_d,\omega,E'_q,U_q,U_d,I_d,\alpha_b^{(1)},\alpha_c^{(1)}\\
\nonumber
&\ &\sim \mathcal{G}\bigg(\alpha_{I}^{(0)}+\dfrac{n}{2}, \\
\nonumber
&\ &\beta_{I}^{(0)}+\dfrac{\sum_{i=1}^{n_3}\big(y_{Iq}-\alpha_b^{(1)}(U_d-E'_d)-\alpha_c^{(1)}(U_q-E'_q)\big)^2+\big(y_{Id}-\alpha_b^{(1)}(U_d-E'_d)
+\alpha_c^{(1)}(U_q-E'_q)\big)^2}{2}\bigg)
\end{eqnarray}
\subsection{Gibbs Sampling Simulation}

\begin{table}[!h]
\caption{Estimated and real value of the parameters for IM model with 1\% measurement error.}
\label{imgibbscompare}
\centering
\begin{tabular}{ c|c|c|c }
\hline
Para.&Real Value &Est. Value (mean)&Error(\%)\\
\hline   
$\alpha_1$&0.0077&0.0076&1.3\\
$\alpha_2$&0.018&0.0188&4.4\\
$\alpha_3$&25&24.28&2.88\\
$\alpha_b$&0.20&0.203&1.5\\
$\alpha_c$&0.80&0.75&6.25\\
\hline
\end{tabular}
\end{table}\par

\begin{table}[!h]
\caption{Estimated and real value of the parameters for IM mode 2\% measurement errorl.}
\label{imgibbscompare2}
\centering
\begin{tabular}{ c|c|c|c }
\hline
Para.&Real Value &Est. Value (mean)&Error(\%)\\
\hline   
$\alpha_1$&0.0077&0.007753&0.69\\
$\alpha_2$&0.018&0.01849&2.7\\
$\alpha_3$&25&24.05&3.8\\
$\alpha_b$&0.20&0.18&10\\
$\alpha_c$&0.80&0.75&6.25\\
\hline
\end{tabular}
\end{table}\par

\begin{table}[!h]
\caption{Estimated and real value of the parameters for IM mode 5\% measurement errorl.}
\label{imgibbscompare3}
\centering
\begin{tabular}{ c|c|c|c }
\hline
Para.&Setting Value &Est. Value (mean)&Error(\%)\\
\hline
$\beta_1$&0.0077&0.007753&0.69\\
$\beta_2$&0.018&0.01849&2.7\\
$\beta_3$&25&24.05&3.8\\
$\alpha_b$&0.20&0.21&5\\
$\alpha_c$&0.80&0.82&2.5\\
\hline
\end{tabular}
\end{table}\par
The simulation results can be found in Fig. \ref{imgibbs}. The real value and estimated value are listed in Table \ref{imgibbscompare}. Compare with the real value, the estimated values are relatively close. The maximum error shown in Table \ref{imgibbscompare} is only 6.25\%. This indicates the proposed GS method works for IM model parameters estimation.
\section{Difficulties \& Conclusion}
Some difficulties during the progress of this topic:\par
1. How to verify the equations used in ZIP and IM model follow normal distribution. This assumption was proven to be correct according to the simulation, but we somehow still need to validate this assumption before we use it.\par
2. Currently we are only estimating the parameters of one ZIP or IM model. It is necessary to consider multi ZIP and MI models connected to the grid. There might be inferences to the estimation. This need to be verify.\par
3. How to identify the composition of ZIP and IM is another interesting topic that can be studied after this report. Most of the works that have been done in composite load modeling need to specify the percentage of ZIP by given the active power load $P_{ZIP}$.\par
Conclusions of this report are:  Bayesian estimation based dynamic load modeling is an executable and reliable composite load modeling approach. This method can overcome some disadvantages in the traditional composite load modeling process. It can provide a robust time-varying dynamic load parameter estimation with higher accuracy and better performance. Future work will be focused on multi composite model identification using hidden Markov Model (HMM), composition of ZIP and IM model, and model clustering using machine learning approach.

\section*{References}

\bibliography{reference}
\section*{Appendix}
The detailed parameters of a 33-bus test feeder can be found as follows:
\begin{table}[!h]
  \centering
   \small
\setlength\tabcolsep{2pt}
  \caption{33-Bus Test Feeder Data}
    \begin{tabular}{ccccccc}
\hline
\hline
  Line &Sending&Receiving&Resistance&Reactance& Real &Reactive\\
Number&     Bus   &      Bus     &    ($\Omega$)       &     ($\Omega$)    &Power&  Power  \\
            &              &                 &                  &                 &(kW)&(kVar)\\
\hline
\hline
    1     & 0     & 1     & 0.0922 & 0.0477 & 100   & 60 \\
    2     & 1     & 2     & 0.493 & 0.2511 & 90    & 40 \\
    3     & 2     & 3     & 0.366 & 0.1864 & 120   & 80 \\
    4     & 3     & 4     & 0.3811 & 0.1941 & 60    & 30 \\
    5     & 4     & 5     & 0.819 & 0.707 & 60    & 20 \\
    6     & 5     & 6     & 0.1872 & 0.6188 & 200   & 100 \\
    7     & 6     & 7     & 1.7114 & 1.2531 & 200   & 100 \\
    8     & 7     & 8     & 1.03  & 0.74  & 60    & 20 \\
    9     & 8     & 9     & 1.04  & 0.74  & 60    & 20 \\
    10    & 9     & 10    & 0.1966 & 0.065 & 45    & 30 \\
    11    & 10    & 11    & 0.3744 & 0.1238 & 60    & 35 \\
    12    & 11    & 12    & 1.468 & 1.155 & 60    & 35 \\
    13    & 12    & 13    & 0.5416 & 0.7129 & 120   & 80 \\
    14    & 13    & 14    & 0.591 & 0.526 & 60    & 10 \\
    15    & 14    & 15    & 0.7463 & 0.545 & 60    & 20 \\
    16    & 15    & 16    & 1.289 & 1.721 & 60    & 20 \\
    17    & 16    & 17    & 0.732 & 0.574 & 90    & 40 \\
    18    & 1     & 18    & 0.164 & 0.1565 & 90    & 40 \\
    19    & 18    & 19    & 1.5042 & 1.3554 & 90    & 40 \\
    20    & 19    & 20    & 0.4095 & 0.4787 & 90    & 40 \\
    21    & 20    & 21    & 0.7089 & 0.9373 & 90    & 40 \\
    22    & 2     & 22    & 0.4512 & 0.3083 & 90    & 50 \\
    23    & 22    & 23    & 0.898 & 0.7091 & 420   & 200 \\
    24    & 23    & 24    & 0.896 & 0.7011 & 420   & 200 \\
    25    & 5     & 25    & 0.203 & 0.1034 & 60    & 25 \\
    26    & 25    & 26    & 0.2842 & 0.1447 & 60    & 25 \\
    27    & 26    & 27    & 1.059 & 0.9337 & 60    & 20 \\
    28    & 27    & 28    & 0.8042 & 0.7006 & 120   & 70 \\
    29    & 28    & 29    & 0.5075 & 0.2585 & 200   & 600 \\
    30    & 29    & 30    & 0.9744 & 0.963 & 150   & 70 \\
    31    & 30    & 31    & 0.3105 & 0.3619 & 210   & 100 \\
    32    & 31    & 32    & 0.341 & 0.5302 & 60    & 40 \\
\hline
\hline
    \end{tabular}%
  \label{33busdata}%
\end{table}%
\end{document}